\documentclass[reqno,12pt]{amsart}
\usepackage{amsfonts,amssymb,amsmath,amscd}
\usepackage{graphicx}
\usepackage{xspace}

\date{Submitted for publication on 16~May 2001,
accepted 15 April 2004.
}

\title[Interlace polynomial]%
{A two-variable interlace polynomial}

\author[Arratia]{Richard Arratia}
\address[Richard Arratia]{Univ.\ of Southern California \\
Department of Mathematics \\
Los Angeles CA 90089-1113, USA}
\email{rarratia@math.usc.edu}

\author[Bollob\'{a}s]{B\'{e}la Bollob\'{a}s$^\dag$}
\address[B\'{e}la Bollob\'{a}s]{Univ.\ of Memphis \\
Department of Mathematical Sciences \\
Memphis TN 38152, USA,
and
Trinity College, Cambridge CB2 1TQ, U.K.
}
\email{bollobas@msci.memphis.edu}

\author[Sorkin]{Gregory B. Sorkin}
\address[Gregory B. Sorkin]{IBM T.J.\ Watson Research Center \\
Department of Mathematical Sciences \\
Yorktown Heights NY 10598, USA}
\email{sorkin@watson.ibm.com}

\thanks{$^\dag$
Supported by NSF grant DMS-9971788.
}


\begin{document}
\bibliographystyle{amsalpha}

\newtheorem{theorem}{Theorem}
\newtheorem{conjecture}[theorem]{Conjecture}
\newtheorem{corollary}[theorem]{Corollary}
\newtheorem{definition}[theorem]{Definition}
\newtheorem{lemma}[theorem]{Lemma}
\newtheorem{proposition}[theorem]{Proposition}
\newtheorem{remark}[theorem]{Remark}

\newcommand{\compl}[1]{#1^c}
\newcommand{\csub}[1]{{\makebox[0in]{$\substack{#1}$}}}
\newcommand{\Gstar}{G^\star}
\newcommand{\Gab}{G^{ab}}
\newcommand{\Ga}{G^a}
\newcommand{\Hab}{H^{ab}}
\newcommand{\uu}{a}
\newcommand{\vv}{b}
\newcommand{\uv}{ab}
\newcommand{\vu}{ba}
\newcommand{\by}{\times}
\newcommand{\qq}{\tilde{q}}
\newcommand{\ilace}{H}
\newcommand{\HH}{H}
\newcommand{\Hset}{\mathcal{H}}
\newcommand{\Pset}{\mathcal{P}}
\newcommand{\Cset}{\mathcal{C}}
\newcommand{\Dset}{\mathcal{D}}
\newcommand{\Dstar}{D^{\star}}
\newcommand{\Cstar}{C^{\star}}
\newcommand{\Dstarp}{{D^{\star}}'}
\newcommand{\Cstarp}{{C^{\star}}'}
\newcommand{\Ab}{\overline{A}}
\newcommand{\ham}{{H}}
\newcommand{\hb}{\overline{H}}
\newcommand{\Ib}{\overline{I}}
\newcommand{\indep}{\operatorname{ind}}
\newcommand{\clique}{\operatorname{cl}}
\newcommand{\num}{\operatorname{num}}
\newcommand{\qh}{q_\HH}
\newcommand{\qhx}{q_\HH(x)}
\newcommand{\R}{\mathbb{R}}
\def\tito{2-in, 2-out }
\newcommand{\ZZ}{\mathbb Z}
\newcommand{\FF}{\mathbb F}
\renewcommand{\Re}{\mathbb R}
\def\cal{\mathcal}
\newcommand{\ckt}{C}
\newcommand{\order}[1]{|{#1}|}
\newcommand{\numeuler}{r_1}
\newcommand{\size}[1]{e({#1})}
\newcommand{\ie}{\textit{i.e.}}
\newcommand{\eg}{\textit{e.g.}}
\newcommand{\restricted}[1]{\: |_{#1}}
\newcommand{\ignorecmd}[1]{}
\newcommand{\eqdef}{\stackrel{\mbox{\scriptsize\rm def}}{=}}
\newcommand{\ind}{\hspace*{1cm}}
\newcommand{\tmbf}[1]{\textbf{\boldmath{{#1}}}}
\newcommand{\isom}{\cong}
\newcommand{\qr}{q_R}
\newcommand{\qn}{q_N}
\newcommand{\sharpp}{\#P}
\newcommand{\vleq}{\preceq}
\newcommand{\vgeq}{\preceq}

\newcommand{\one}{{\mathbf{1}}}
\newcommand{\zero}{{\mathbf{0}}}
\newcommand{\bs}{\backslash}
\newcommand{\union}{\cup}
\newcommand{\odd}{\operatorname{odd}}
\newcommand{\nullity}{vertex-nullity\xspace}
\newcommand{\rank}{vertex-rank\xspace}
\newcommand{\nbd}{\Gamma}

\begin{abstract}
We introduce a new graph polynomial in two variables.  This ``interlace'' polynomial can be computed in two very different ways.  The first is an expansion analogous to the state space expansion of the Tutte polynomial; the significant differences are that our expansion is over vertex rather than edge subsets, and the rank and nullity employed are those of an adjacency matrix rather than an incidence matrix.

The second computation is by a three-term reduction formula involving a graph pivot; the pivot arose previously in the study of interlacement and Euler circuits in four-regular graphs.

We consider a few properties and specializations of the two-variable interlace polynomial.  One specialization, the ``\nullity interlace polynomial'', is the single-variable interlace graph polynomial we studied previously, closely related to the Tutte-Martin polynomial on isotropic systems previously considered by Bouchet.  Another, the ``\rank interlace polynomial'', is equally interesting.  Yet another specialization of the two-variable polynomial is the independent-set polynomial.
\end{abstract}

\maketitle

\section{The interlace polynomial}
\label{intro}

In~\cite{ABS00-AMSformat,abs-jctb},
we introduced a single-variable ``interlace'' graph polynomial.
It emerged that the interlace polynomial
could be regarded as
a special case of the Tutte-Martin polynomial of an isotropic system,
as discussed briefly in Section~\ref{sec:special} here
and more fully in~\cite{abs-jctb}.
We defined the polynomial by a recurrence relation, and
Balister, Bollob\'as, Cutler and Pebody~\cite{BBCP02}
used a property of the nullities of the adjacency matrices of
the graphs in question to resolve a conjecture posed in~\cite{ABS00-AMSformat}.
The linear-algebraic approach of~\cite{BBCP02},
extended to embrace the matrix ranks as well as nullities,
led us to the two-variable polynomial introduced here.
The two-variable interlace polynomial appears to be something entirely new;
in particular, we are aware of no two-variable versions of
the Tutte-Martin polynomial or other closely related polynomials.
The interlace polynomial is an entirely different object
from the Tutte polynomial,
but the two have an extremely interesting structural similarity
which immediately suggests a family of additional polynomials
for further exploration.
We will now define the two-variable interlace polynomial.

Given a graph $G$ with vertex set $V(G)$,
for any subset $S \subseteq V(G)$, let $G[S]$ be
the subgraph of $G$ induced by $S$.
We allow graphs with loops on their vertices,
and we also allow the {\em null graph} with no vertices,
but we do not allow multiple loops or multiple edges.
Write $\mathcal{G}$ for the set of graphs including the null graph.

For a matrix $A$ over $\FF_2$,
let $n(A)$ be the nullity of $A$ and $r(A)$ its rank.
Abusing notation slightly, for a graph $G$,
$n(G)$ and $r(G)$ will denote the nullity and rank of its
adjacency matrix,
so $n(G)+r(G) = |V(G)|$.
(For example,
if $S$ is the empty set then $G[S]$ is the null graph of rank and
nullity~$0$.)
We remark that for loopless graphs~$G$, $r(G)$ is always even,
as the rank of a zero-diagonal symmetric matrix.

We define the two-variable interlace polynomial
$q(G;x,y)$ of a
graph $G$ of order $n$ as a sum of $2^n$ terms:
\begin{align}
q(G;x,y) &= \sum_{S \subseteq V(G)} (x-1)^{r(G[S])} (y-1)^{n(G[S])} ,
\label{explicit}
\end{align}
the sum taken over all subsets including $S=\emptyset$ and $S=V(G)$.
For convenience we define the monomial
\begin{align*}
m(H) &= (x-1)^{r(H)} (y-1)^{n(H)} ,
\end{align*}
so that
\begin{align}
q(G) &= \sum_S m(G[S]) .
\label{explicitm}
\end{align}

This ``state space'' expansion of
the two-variable interlace polynomial may be seen as an analogue of the Tutte
polynomial given as
\begin{align}
T(G;x,y) &= \sum_{F \subseteq E} (x-1)^{r(E)-r(F)} (y-1)^{n(F)} .
\label{Tutte}
\end{align}
The significant differences are that our sum
is over induced subgraphs of $G$ (as given by vertex subsets) rather than
arbitrary subgraphs on the full vertex set (as given by edge subsets),
and our rank and nullity are those of the
subgraph's adjacency matrix rather than its incidence matrix.
(In the context of the Tutte polynomial,
a graph's rank is normally defined
as the number of vertices minus the number of components,
and it is easy to check
that this is the $\FF_2$-rank of the incidence matrix.)
That our subgraph rank appears positively rather than subtracted
from the rank of the whole is not significant, as it can be adjusted
by a change of variables.  That is, if one prefers the polynomial
$\sum_{S \subseteq V(G)} (x-1)^{r(G)- r(G[S])} (y-1)^{n(G[S])}$,
it is just $(x-1)^{r(G)} q(\frac{x}{x-1},y)$.

A surprising basic property of this polynomial is that
for loopless graphs it satisfies a
three-term \emph{reduction formula}, as per Theorem~\ref{reduction},
and for looped graphs, a pair of reductions, per Theorem~\ref{thm:red12};
we prove these results in Section~\ref{sec:reduction}.
In Section~\ref{sec:properties} we show a pair of properties of the
two-variable interlace polynomial.
In Section~\ref{sec:special} we describe the polynomial's specializations
to the vertex-rank polynomial and the vertex-nullity polynomial,
and in Section~\ref{sec:independent},
its specialization to the independent-set polynomial.
We calculate the polynomial on some basic graphs in Section~\ref{sec:graphs}.
We conclude in Sections~\ref{otherpolys} and \ref{open}
with thoughts on generalizations of the polynomial,
some of which seem likely to prove interesting,
and with some open problems.

We now proceed to the interlace polynomial's reduction formula.

\section{Reduction formula}
\label{sec:reduction}

We begin
by showing a reduction
involving a ``pivoting'' operation on an edge $ab$ of $G$
for which both $a$ and $b$ are loopless.
It is natural to start with this case because it provides a
recursive definition of $q(G)$ for the natural class of loopless graphs~$G$.
In the following subsection, we show how to
reduce on a looped vertex $a$ of~$G$,
completing a recursive definition of~$q(G)$ for arbitrary graphs.

\subsection{The pivot and reduction}
\label{pivoting}

As in~\cite{abs00,abs-jctb,Bouchet-mm3},
and related to Kotzig's transformations on Euler tours~\cite{Kotzig69},
for a graph $G$ and an ordered pair
$ab=(a,b)$ of distinct vertices of $G$, we define the
{\em pivot operation on $ab$} mapping $G$ into $ G^{ab}$
as follows. We say that
two vertices $x, y$ of $G$ are {\em distinguished} by $\{a, b\}$
if $x, y \notin \{a, b\}$ and $x, y$ have distinct non-empty
neighborhoods in $\{a, b\}$. Let
$G^{ab}$ be the graph with vertex set $V(G)$ in which $xy$ is an edge
if either $xy\notin E(G)$ and $x$ and $y$ are distinguished by $\{a,
b\}$, or else $xy\in E(G)$ and $x$ and $y$ are not distinguished by
$\{a, b\}$.

Let us spell out this definition in detail.
Partition the vertices other than $\uu$ and $\vv$ into four classes:
\begin{description}
\item[C1] vertices adjacent to both $\uu$ and $\vv$;
\item[C2] vertices adjacent to $\uu$ alone;
\item[C3] vertices adjacent to $\vv$ alone;
and
\item[C4] vertices adjacent to neither $\uu$ nor $\vv$.
\end{description}

\begin{definition}[Pivot]
A graph $G$ is pivoted on vertices $a,b$ to obtain $G^{ab}$ as follows.
For any vertex pair $xy$ where $x$ is in one of the classes~(C1)--(C3)
and $y$ is in a different class~(C1)--(C3),
the pair $xy$ is ``toggled'':
if it is an edge of $G$ it is not an edge of $G^{ab}$,
and if it is not an edge of $G$ then it is an edge of $G^{ab}$.
All other pairs of vertices are adjacent in $G^{ab}$ iff they are adjacent
in~$G$.
\end{definition}

Trivially, pivoting and restriction to a subgraph satisfy a
commutative law: for any $S \subseteq V(G)$ with $a,b \in S$,
\begin{align}
G^{ab}[S] &= (G[S])^{ab} .
\label{commute}
\end{align}

Although pivoting is defined for any vertex pair $ab$,
we shall only exploit it in cases where
$ab$ is an edge and $a$ and $b$ are both loopless.

We shall write the adjacency matrices of $G$ and
$G^{ab}$ with rows and columns put
into six groups according to their relations to $a$ and $b$. The first
group consists of $a$ alone, and the second of $b$ alone; groups three
to six are the four classes above.
Write $\one$ for an all-1 row or column
vector of whatever dimension and likewise $\zero$ for an all-0 vector.
Then, for a graph $G$ with loopless vertices $a$ and $b$ and
having an edge $ab$,
the adjacency matrix of $G$ is of the form
\begin{align}
A(G) &=
\left( \begin{array}{cccccc}
0&1&\one&\one&\zero&\zero \\
1&0&\one&\zero&\one&\zero \\
\one&\one&M_{11}&M_{12}&M_{13}&M_{14} \\
\one&\zero&M_{21}&M_{22}&M_{23}&M_{24} \\
\zero&\one&M_{31}&M_{32}&M_{33}&M_{34} \\
\zero&\zero&M_{41}&M_{42}&M_{43}&M_{44} ,
\end{array} \right)
\label{adjmat}
\\
\intertext{%
where in all cases $M_{ji}$ is the transpose of $M_{ij}$.
Then for such a graph the adjacency matrix of $G^{ab}$ is
}
A(\Gab) &=
\left( \begin{array}{cccccc}
0&1&\one&\one&\zero&\zero \\
1&0&\one&\zero&\one&\zero \\
\one&\one&M_{11}&\compl{M_{12}}&\compl{M_{13}}&M_{14} \\
\one&\zero&\compl{M_{21}}&M_{22}&\compl{M_{23}}&M_{24} \\
\zero&\one&\compl{M_{31}}&\compl{M_{32}}&M_{33}&M_{34} \\
\zero&\zero&M_{41}&M_{42}&M_{43}&M_{44}
\end{array} \right)
\label{adjmatab}
\end{align}
where  $M^c$ denotes the complement of the matrix~$M$.

\begin{lemma}
\label{essence}
For any graph $G$ with an edge $ab$, with $a$ and $b$ both loopless,
$r(G-a) = r(\Gab-a)$ and $r(G) = r(\Gab-a-b)+2$;
equivalently,
$n(G-a) = n(\Gab-a)$ and $n(G) = n(\Gab-a-b)$.
\end{lemma}

\begin{proof}
With the adjacency matrix $A=A(G)$ represented as in~\eqref{adjmat},
add row 2 to each row in the 3rd and 4th groups.
Repeat the same operations on columns to give a matrix
\begin{align*}
A' &=
\left( \begin{array}{cccccc}
0&1&\zero&\zero&\zero&\zero \\
1&0&\one&\zero&\one&\zero \\
\zero&\one&M_{11}&\compl{M_{12}}&\compl{M_{13}}&M_{14} \\
\zero&\zero&\compl{M_{21}}&M_{22}&\compl{M_{23}}&M_{24} \\
\zero&\one&\compl{M_{31}}&\compl{M_{32}}&M_{33}&M_{34} \\
\zero&\zero&M_{41}&M_{42}&M_{43}&M_{44}
\end{array} \right) .
\end{align*}
Note that, outside of the neighborhoods of $a$ and $b$
(the first two rows and columns),
$A'$ is the adjacency matrix of $\Gab$; moreover,
since these linear operations are invertible, $r(A)=r(A')$.

To prove the first assertion, discard the first row and column of
$A$ to yield $A \bs a$ and similarly that of $A'$ to obtain $A' \bs a$.
Since they did not use row or column~1 (vertex $a$), the same linear
transformations as before map $A \bs a$ to $A' \bs a$,
showing that $r(A \bs a) = r(A' \bs a)$.
But $A \bs a = A(G-a)$
and $A' \bs a = A(\Gab-a)$, therefore $r(G-a) = r(\Gab-a)$.

To prove the second assertion, further transform $A'$ by adding
row 1 to each row in the 3rd and 5th groups, and repeating for columns,
to obtain
\begin{align*}
A'' &=
\left( \begin{array}{cccccc}
0&1&\zero&\zero&\zero&\zero \\
1&0&\zero&\zero&\zero&\zero \\
\zero&\zero&M_{11}&\compl{M_{12}}&\compl{M_{13}}&M_{14} \\
\zero&\zero&\compl{M_{21}}&M_{22}&\compl{M_{23}}&M_{24} \\
\zero&\zero&\compl{M_{31}}&\compl{M_{32}}&M_{33}&M_{34} \\
\zero&\zero&M_{41}&M_{42}&M_{43}&M_{44}
\end{array} \right) .
\end{align*}
Because these linear transformations are all invertible, they preserve
rank and nullity: $r(A) = r(A'')$.
Referring to \eqref{adjmatab},
note that $A'' \bs \{a,b\}$ is the adjacency matrix of $\Gab-a-b$.
The first and second rows in $A''$ are linearly independent of one another
and of all other rows, so deleting them
reduces the rank by~2.
After deletion of these rows the first two columns are all-zero, so
deleting them
does not change the rank;
it follows that $r(A'') = r(A'' \bs \{a,b\}) +2$.
That is, $r(G) = r(A) = r(A'') = r(\Gab-a-b)+2$.
\end{proof}

Note that if either $a$ or $b$ has a loop,
the top-left submatrix of $A(G)$ differs from that in~\eqref{adjmat},
resulting in a different ``border'' in the matrix $A''$, so that
the border's deletion changes the nullity unpredictably.
That is, if there is a loop at either $a$ or $b$,
$n(G)-n(\Gab-a-b)$ may be 0 rather than~2.
Thus, the looped case is dealt with in Section~\ref{complementing}.

\begin{theorem}
\label{reduction}
For any edge $ab$ of a graph $G$, where $a$ and $b$ are both loopless,
\begin{align*}
q(G) = q(G-a) + q(\Gab-b) + ((x-1)^2-1) q(\Gab-a-b) .
\end{align*}
\end{theorem}

\begin{proof}
For $S$ ranging over subsets of $V(G) \bs \{a,b\}$,
by~\eqref{explicitm},
\begin{align}
q(G) &=
\sum_S \{ m(G[S]) + m(G[S \union a])
 + m(G[S \union b]) + m(G[S \union \{a,b\}]) \}
\label{e1}
\end{align}
while
\begin{align}
\lefteqn{
q(G-a)+q(\Gab-b)+((x-1)^2-1) q(\Gab-a-b)
} &
\notag \\ &=
 \sum_S \{ m((G-a)[S]) + m((G-a)[S \union b]) \}
\notag \\ & \phantom{=}
 + \sum_S \{ m(\Gab-b)[S]) + m((\Gab-b)[S \union a]) \}
\notag \\ & \phantom{=}
 + \sum_S ((x-1)^2-1) m((\Gab-a-b)[S]) .
\label{e2}
\end{align}
To show that \eqref{e1} and \eqref{e2} are equal,
we will show equality of their terms for each~$S$.
Two terms of \eqref{e1} directly match their counterparts in \eqref{e2}:
$m(G[S]) = m((G-a)[S])$ and
$m(G[S \union b]) = m((G-a)[S \union b])$.
A third equality follows from the first part of Lemma~\ref{essence}:
\begin{align*}
m(G[S \union a])
 &= m(G[S \union \{a,b\}]-b)
 \\& = m((G[S \union \{a,b\}])^{ab} -b)
 \\& = m(\Gab[S \union a])   \qquad \text{(by \eqref{commute})}
 \\& = m((\Gab-b)[S \union a]) .
\end{align*}
The final equality, between a single term from \eqref{e1}
and two terms from \eqref{e2},
follows from the second part of Lemma~\ref{essence}:
\begin{align*}
m(G[S \union \{a,b\}])
 &= (x-1)^2 m((G[S \union \{a,b\}])^{ab}-a-b)
 \\& = (x-1)^2 m(\Gab[S])
 \\& = m(\Gab[S]) + ((x-1)^2-1) m(\Gab[S])
 \\& = m((\Gab-b)[S]) + ((x-1)^2-1) m((\Gab-a-b)[S]) .
\end{align*}
\end{proof}

\subsection{Local complementation}
\label{complementing}

In the case graphs with loops,
there may not always be an edge $ab$, with both $a$ and $b$ loopless,
to which Theorem~\ref{reduction} (derived from Lemma~\ref{essence})
may be applied.
In this case, though, there must be a looped vertex,
and we may apply a different reduction instead.

Bouchet~\cite{Bouchet-mm3} defines the
``local complement'' $\Ga$ of graph $G$ on vertex $a$
by complementing
(toggling the presence or absence of all edges, including loops)
the subgraph of $G$ induced by the neighborhood of $a$, while
keeping the graph otherwise unchanged.
To be precise, let $\nbd(a)$ be the set of neighbors of $a$;
in particular $a \in \nbd(a)$ iff there is a loop on~$a$.

\begin{definition}[Local complementation]
A graph $G$ is locally complemented on a vertex $a$ to yield $\Ga$,
where $\Ga$ is equal to $G$ except that
$\Ga[\nbd(a)] = \overline{G[\nbd(a)]}$.
\end{definition}

In notation like that of \eqref{adjmat}
and \eqref{adjmatab} but only distinguishing
a vertex~$a$ having a loop,
its neighbors, and its non-neighbors,
we may write
\begin{align}
A(G) &=
\left( \begin{array}{ccc}
1&\one&\zero\\
\one&M_{11}&M_{12} \\
\zero&M_{21}&M_{22} \\
\end{array} \right)
\label{adjmat2}
\\ \intertext{and}
A(\Ga) &=
\left( \begin{array}{ccc}
1&\one&\zero\\
\one&\compl{M_{11}}&M_{12} \\
\zero&M_{21}&M_{22} \\
\end{array} \right) .
\label{adjmat2a}
\end{align}

Incidentally, it is observed in~\cite{Bouchet-mm3} that a pivot is equal to a
composition of local complementations, $\Gab = ((\Ga)^b)^a$,
followed by a swap of the labels $a$ and $b$.
(This is correct as stated for our version of pivoting,
which differs from Bouchet's by a label swap.)

\begin{lemma}
\label{essence2}
For any graph $G$ with a looped vertex~$a$,
$r(G) = r(\Ga-a)+1$,
and equivalently $n(G) = n(\Ga-a)$.
\end{lemma}

\begin{proof}
With the adjacency matrix $A=A(G)$ represented as in~\eqref{adjmat2},
let $A'$ be obtained by adding row~1 of $A(G)$ to each row in the
second class, and repeating for columns; thus
\begin{align*}
A' &=
\left( \begin{array}{ccc}
1&\zero&\zero\\
\zero& \compl{M_{11}} &M_{12} \\
\zero&M_{21}&M_{22} \\
\end{array} \right) .
\end{align*}
Referring to \eqref{adjmat2a},
note that $A(\Ga-a) = A(\Ga) \bs a = A' \bs a$.
The linear transformations are invertible, so $r(A')=r(A)$.
The first row of $A'$ is independent of the others, so
deleting it decreases the rank by~1,
and what remains of the first column is all-zero,
so deleting it does not change the rank;
thus $r(A' \bs a) = r(A')-1$.
We conclude that $r(\Ga-a) = r(A' \bs a) = r(A')-1 = r(A)-1 = r(G)-1$.
\end{proof}

\begin{theorem}
\label{thm:red12}
For a graph $G$,
for any edge $ab$ where neither $a$ nor $b$ has a loop,
\begin{align}
q(G) &= q(G-a) + q(\Gab-b) + ((x-1)^2-1) q(\Gab-a-b) ,
 \label{red1}
\\ \intertext{and for any looped vertex $a$,}
q(G) &= q(G-a) + (x-1) q(\Ga-a) ,
 \label{red2}
\end{align}
\end{theorem}

\begin{proof}
Equation~\eqref{red1} is precisely Theorem~\ref{reduction},
repeated here to give one comprehensive theorem.

To prove \eqref{red2},
as in the proof of Theorem~\ref{reduction}, we show that for
each $S \subseteq V(G) \bs a$, we have equality between the
corresponding summands of
\begin{align*}
q(G) &= \sum_S \{ m(G[S]) + m(G[S \union a]) \}
\\ \intertext{and those of}
q(G-a) + (x-1) m(\Ga-a) &= \sum_S m((G-a)[S]) + (x-1) \sum_S m((\Ga-a)[S]) .
\end{align*}
The first terms, $m(G[S])$ and $m((G-a)[S])$, are identical.
By Lemma~\ref{essence2},
$m(G[S \union a]) = (x-1) \: m((G[S \union a])^a -a)
 = (x-1) \: m((\Ga-a)[S])$,
which completes the proof.
\end{proof}

The reduction formulas give an alternative characterization of the
two-variable interlace polynomial.
We write $E_n$ for the empty graph on $n$ vertices.

\begin{corollary}
\label{implicit}
The two-variable interlace polynomial
defined by \eqref{explicit}
is the unique map
$q: \mathcal{G} \to {\mathbb Z}[x,y]$ that satisfies the reduction
formulas (\ref{red1},\ref{red2}) and the boundary conditions
$q(E_n)= y^n$, $n=0, 1, \ldots$.
\end{corollary}

\begin{proof}
By Theorem~\ref{thm:red12}, the two-variable interlace polynomial $q(G)$
defined by \eqref{explicit} satisfies the reduction formulas
(\ref{red1},\ref{red2})
and from \eqref{explicit} it is immediate that it also satisfies
the boundary conditions.
Uniqueness follows from (\ref{red1},\ref{red2})
by induction on the order of the graph.
\end{proof}

\section{Two properties}
\label{sec:properties}

As will be shown in the next section,
the single-variable interlace polynomial defined in~\cite{abs00},
which there was denoted as $q(G;x)$,
is a special case of the present two-variable interlace polynomial,
and we will denote it here as $\qn(G;y)$.
(The reason for this notation, and the definition of $\qn$,
will be given in the next section.)

Since \cite{abs00} showed that
this single-variable interlace polynomial satisfies the identity
$\qn(G;y) = \qn(\Gab;y)$
(at least for loopless graphs, all that \cite{abs00} considered),
something of the same sort might be expected for the
two-variable polynomial.
In fact it is not generally true that $q(G;x,y)$ is equal to $q(\Gab;x,y)$
--- a counterexample is the path of length~3, pivoted on the middle edge
--- but instead we have the following proposition.

\begin{proposition}
For any graph $G$ with edge $ab$,
\begin{align*}
q(G-a)-q(G-a-b) = q(\Gab-a)-q(\Gab-a-b) .
\end{align*}
\end{proposition}

\begin{proof}
By the interlace polynomial's definition,
with sums taken over subsets $S \subseteq V(G) \bs \{a,b\}$,
\begin{align*}
q(G-a)-q(G-a-b) &= \sum_S (m(G[S]) + m(G[S \union b]))- \sum_S m(G[S])
\\ &= \sum_S m(G[S \union b])
\\ &= \sum_S m(\Gab[S \union b])
\end{align*}
by Lemma~\ref{essence}.  From this point symmetry completes the proof.
\end{proof}

As with the earlier single-variable polynomial $\qn(G;y)$,
the present two-variable polynomial obeys a simple product rule.
For graphs $G_1 = (V_1, E_1)$ and $G_2 = (V_2, E_2)$ with
disjoint vertex sets, $V_1 \cap V_2 = \emptyset$,
let $G_1 \cup G_2$ denote their disjoint union,
$(V_1 \cup V_2, E_1 \cup E_2)$.

\begin{proposition}
If $G_1$ and $G_2$ are graphs on disjoint vertex sets,
and $G_1 \cup G_2$ their disjoint union,
then
$q(G_1 \cup G_2) = q(G_1) q(G_2)$.
\end{proposition}

\begin{proof}
As in~\cite{abs00}, the relation can be proved through the reduction,
taking advantage of the fact that pivots in one component do not affect
the other.
Alternatively, it is easy to see that
the sum~\eqref{explicit} is the product of the corresponding
sums over all $S_1 \subseteq V(G_1)$ and $S_2 \subseteq V(G_2)$.
\end{proof}

\section{Specializations of the interlace polynomial}
\label{sec:special}

Specializing $q(G;x,y)$ by setting $x=2$
(or, for loopless graphs, where the rank is always even, $x=0$)
causes the ``rank'' term to disappear from
\eqref{explicit}, and so we call this polynomial
the ``\nullity interlace polynomial'',
\begin{align}
\qn(G;y) &= q(G;2,y)
= \sum_{S \subseteq V(G)} (y-1)^{n(G[s])} .
\label{qnexpansion}
\end{align}

As alluded to in the previous section, this is
the single-variable interlace polynomial studied in~\cite{abs00}.
To see this, note that
substituting $x=2$ into \eqref{red1} gives
the reduction $\qn(G;y) = \qn(G-a;y)+\qn(\Gab-b;y)$, precisely the
reduction that defined the
single-variable polynomial in~\cite{abs00}
(which was defined only for loopless graphs).
The boundary conditions $\qn(E_n;y)=y^n$ also match,
concluding the equivalence.

\begin{corollary}
The single-variable interlace polynomial of~\cite{abs00},
defined there by a reduction formula,
has the explicit expansion~\eqref{qnexpansion}.
\end{corollary}

This single-variable interlace polynomial
$\qn(G)$ is related to
the Martin polynomial and circuit partition polynomials
\cite{M2,Bouchet87-isotropic-systems,
Bouchet88-graphic-presentations,Bouchet91,
E1,JoEM00,bollobas-circuit,abs-jctb}.
Bouchet \cite{Bouchet-personal-comm00}
recognized that
the single-variable interlace graph polynomial of \cite{abs00}
was a specialization of
the Tutte-Martin polynomial of an isotropic system
(introduced by Bouchet
in \cite{Bouchet87-isotropic-systems}
and generalizing the Martin polynomial),
and this connection was clarified and made explicit
by Aigner and van der Holst in~\cite{AvdH}.
\cite{AvdH},~written after the present work's submission,
also proved a conjecture
from \cite{ABS00-AMSformat} that $\qn(G;-1)$ is always of the form $\pm 2^s$
(previously proved in \cite{BBCP02}),
independently derived
the expansion \eqref{qnexpansion} of the nullity polynomial,
and introduced a related polynomial.
However, in contrast to the Tutte-like two-variable
graph polynomial of the present paper,
we are not aware of any two-variable generalization
of the Tutte-Martin polynomial.

Analogously to the nullity polynomial $\qn(G;y)$,
there is also a single-variable ``\rank interlace polynomial''
\begin{align*}
\qr(G;x) &= q(G;x,2) = \sum_{S \subseteq V(G)} (x-1)^{r(G[S])} ,
\end{align*}
and it appears to be equally interesting.
First, the \rank polynomial distinguishes graphs of small order better
than the \nullity polynomial.
For example, the rank polynomial distinguishes all 11 simple graphs of order~4,
where the nullity polynomial takes on only 8 distinct values.
(And the rank polynomial distinguishes all 90 looped graphs of order~4,
where the nullity polynomial takes on only 17 distinct values).
At order 5 there are 34 non-isomorphic simple graphs:
the rank polynomial takes 33 values, and
the nullity polynomial only~17.
(There are 544 non-isomorphic looped graphs,
the rank polynomial takes 541 values, and
the nullity polynomial only 41.)
Similarly for trees: the nullity polynomial fails to distinguish one pair of
trees of order 8 and two pairs of order~9;
the rank polynomial distinguishes all trees of orders 8 and~9.

In~\cite{abs-jctb} it was shown that certain basic graph parameters could be
read out from the \nullity polynomial, namely the order, the component number,
the edge-independence number, and an upper bound on the (vertex) independence
number.
The \rank polynomial, too, gives the order.

\begin{remark}
For any graph $G$ of order $n$, $\qr(G;2) = 2^n$.
\end{remark}
\begin{proof}
The formula in \eqref{explicit} reduces to a sum, over
all $2^n$ subgraphs of $G$, of $1$ raised to a power.
\end{proof}

As per the following proposition,
the maximum degree of either variable in the two-variable interlace
polynomial is unchanged by ``removing'' the other variable (substituting~2).
\begin{proposition}
For any graph $G$, $\deg_x(q(G;x,y)) = \deg(\qr(G;x))$, and
$\deg_y(q(G;x,y)) = \deg(\qn(G;y))$,
where $\deg_x$ (respectively $\deg_y$) denotes the maximum degree of $x$
(resp.~$y$) in the polynomial.
\end{proposition}

\begin{proof}
We will prove the statement for the \rank polynomial;
that for the \nullity polynomial is proved identically.
Since $\qr(G;x) = q(G;x,2)$,
$\deg(\qr(G;x)) \leq \deg_x(q(G;x,y))$.
Consider any $S \subseteq V(G)$ contributing to \eqref{explicit} a term of
the maximum $x$-degree, degree~$k$.  But
$\qr(G;x) = \sum_{S \subseteq V(G)} (x-1)^{r(G[S])}$,
so each such term here also has $x$-degree~$k$.
For each such~$S$ the coefficient of $x^k$ is 1:
there is no cancellation, and so the $x$-degree is $k$ in $\qr$
as it was in $q$.
\end{proof}

In \cite{abs-jctb} we showed that $\deg(\qn(G)) \geq \indep(G)$, that is,
$\deg(\qn(G))$ is an upper bound on the independence number.
(While it was proved there only for loopless graphs,
the same result for looped graphs follows from \eqref{indeppoly},
in the next section.)
We showed graphs for which the absolute gap
$\deg(\qn(G))-\indep(G)$ was arbitrarily large,
but it remained open whether the ratio $\deg(\qn(G)) / \indep(G)$ could be
made arbitrarily large.  Indeed it can.

Let $\ham_d$ be the $d$-dimensional Hamming cube (so $|\ham_d|=2^d$)
and let $\hb_d$ be its complement.
Then we have the following.

\begin{remark}
The Hamming cube and its complement have (respectively) clique and
independence numbers $\clique(\ham_d) = \indep(\hb_d) = 2$.
For $d>1$ odd, their nullities are $n(\ham_d)=0$ and $n(\hb_d) = 2^d/2$,
and for $d$ even, $n(\ham_d) = 2^d/2$ and $n(\hb_d) = 0$.
\end{remark}

For $\hb_d$ with $d>1$ odd, the set $S=V(\hb_d)$ is a witness that
$\deg(\qn(\hb_d)) \geq 2^d/2$, and therefore
$\deg(\qn(\hb_d))/\indep(\hb_d) \geq 2^d/4$.

\begin{proof}
That $\indep(\hb_d)=\clique(\ham_d)=2$ is immediate from
the structure of the Hamming cube.
We next derive the nullity of $\ham_d$, from which that of $\hb_d$
follows quickly.
Let the adjacency matrix of $\ham_d$ be $A_d$ and that of
$\hb_d$ be $\Ab_d$.

We first show that $A_d$ is self-inverse, for odd~$d$.
The dot product of the $i$th and $j$th rows of $A_d$ is the parity of
the number of vertices at Hamming distance 1 to both $i$ and $j$ in~$\ham_d$.
If $i=j$, then there are $d$ such vertices, for parity~1.
If the distance between $i$ and $j$ is more than~2, then there are no
such vertices: parity~0.
And if $i$ and $j$ are at distance exactly~$2$,
then the vertices in question must agree with both $i$ and $j$ where
those two agree, agree with $i$ in one of the two coordinates where
it differs from~$j$, and agree with $j$ in the other;
there are 2 such vertices, for parity~0.
When $d$ is even everything is the same except in the case $i=j$,
where the $d$ common neighbors mean parity~0, and thus $A_d^2 = \zero$.

The situation is similar for $\Ab_d$, with the parities reversed:
Modulo~2, $\Ab_d = A_d + \one + I$,
where $\one$ is the all-1 matrix.
Then for $d$ odd $\Ab_d^2 = I$, and for $d$ even $\Ab_d^2 = \zero$.

For $d$ odd, invertibility of $A_d$ means it is of full rank.
The self-invertibility of $A_d$ for $d$ odd can also be used to
show that $A_d$ is of half the full rank for $d$ even:
``Gluing'' together two copies of $\ham_{d-1}$ to make an $\ham_d$,
$A_d =
 \left( \begin{array}{cc} A_{d-1} & I \\ I & A_{d-1} \end{array} \right)$,
self-invertibility of $A_{d-1}$ means the second
set of rows is simply $A_{d-1}$ times the first set.
Since the second set is a linear combination of the first set,
the rank of the whole matrix is at most $2^{d-1}$;
the presence of a block $I$ means this rank is achieved.

Similarly, for $\hb_d$ we have
$\Ab_d = \left(
  \begin{array}{cc} \Ab_{d-1} & \Ib \\ \Ib & \Ab_{d-1} \end{array} \right)$,
where $\Ib = \one+I$.
Here we find that for $d$ odd,
the second set of rows is just
$\Ab_{d-1} \Ib = \Ib \; \Ab_{d-1}$
times the first set,
and conclude that $r(\Ab_d) = 2^{d-1}$.
\end{proof}

\section{Counting independent sets}
\label{sec:independent}

There are other interesting specializations of the 2-variable interlace
polynomial.
Evaluating at $y=1$ means that $(y-1)^{n(G[H])}=0$ except
when $n(G[H])=0$ giving $(y-1)^{n(G[H])}=1$.
In particular, then,
$q(G;2,1)
 = \sum_{H \subseteq G} (2-1)^{r(H)} (1-1)^{n(H)}
 = \sum_{H: \; n(H)=0} 1^{r(H)}$
counts full-rank induced subgraphs of~$G$.

Similarly, $q(G;1,2) = \sum_{H: \; r(H)=0} 1^{n(H)}$
counts the independent sets of~$G$
(including the empty set),
a problem that has received widespread attention.
In particular, it is known that counting independent sets
(computing the independence number) is
\sharpp-complete even for low-degree graphs~\cite{DyerGreenhill00},
so it follows that it is \sharpp-hard to compute the
two-variable interlace polynomial
(in particular at the point $(x,y)=(1,2)$)
and the ``rank'' interlace polynomial (at $x=1$).
In fact, it is hard to compute the independence number even
approximately~\cite{DFJ98},
and so the interlace polynomial must be hard even to approximate
at these points.

Given the similarity to the Tutte polynomial, which is hard to compute
almost everywhere (\cite{Jaeger90}, see also \cite{Welsh93} for a survey),
and given the variety of structures counted by the interlace polynomial
(see~\cite{abs-jctb}), with counting typically being \sharpp-hard,
it is anything but surprising that the interlace polynomial
is computationally hard.
However, it was a question left unresolved in~\cite{abs-jctb}, and in fact
we still do not have a proof that computing the ``nullity'' polynomial is
\sharpp-hard.
Moreover, in analogy with the Tutte polynomial, it would be of interest
to show that the interlace polynomial is hard to compute at almost all
points $(x,y)$.

A particularly interesting evaluation is
\begin{align}
q(G;1,1+\lambda)
 = \sum_{H \subseteq G} 0^{r(H)} \lambda^{n(H)}
 = \sum_{H: \; r(H)=0} \lambda^{|H|}
 = I(G;\lambda) ,
\label{indeppoly}
\end{align}
the independent-set polynomial
(the sum over all $k \geq 0$ of $\lambda^k$ times the number of
independent sets of cardinality $k$).
It is well known that further quantities of interest can be computed
from the independent-set polynomial and its derivatives.
For example from \eqref{indeppoly} it is clear that
$\frac{\partial}{\partial \lambda} q(G;1,1+\lambda)
 = \sum_{H: \; r(H)=0} |H| \lambda^{|H|-1}$.
For $\lambda=1$ this is just $\sum_{H: \; r(H)=0} |H|$,
so that
$\frac{\partial q(G)}{\partial y} (1,2)$
is the sum of the sizes of all independent sets.

\section{Polynomials of some basic graphs}
\label{sec:graphs}

We compute the interlace polynomial of some basic graphs, notably
complete graphs $K_n$, complete bipartite graphs $K_{m,n}$,
and paths $P_n$ of length~$n$.
\begin{proposition}
For all $n$ and $m$ we have
\allowdisplaybreaks{
\begin{align*}
q(E_n) &= y^n
\\
q(K_n) &= \frac12 \left( x^n+(2-x)^n \right) +
          \frac12 \left( \frac{y-1}{x-1} \right) \left( x^n-(2-x)^n \right)
\\
q(K_{m,n})
 &= \frac{(x-1)^2}{(y-1)^2} \left( (y^m-1)(y^n-1) \right)
      + y^m + y^n -1
\\
q(P_n)
 &= \frac12 \left(y+ \frac{3y+2x(x-2)}{\sqrt{1+4(y+x(x-2))}} \right)
    \left(\frac{1+ \sqrt{1+4(y+x(x-2))}}{2} \right)^n
   \\ & \phantom{=}
  + \frac12 \left(y- \frac{3y+2x(x-2)}{\sqrt{1+4(y+x(x-2))}} \right)
    \left(\frac{1- \sqrt{1+4(y+x(x-2))}}{2} \right)^n
\end{align*}
}
\end{proposition}

\begin{proof}
That $q(E_n)=y^n$ is immediate from \eqref{explicit} and also figured into the
boundary condition in Corollary~\ref{implicit}.

For $K_n$ we have
\begin{align*}
q(K_n) &= \sum_{k \leq n} \binom{n}{k} (y-1)^{n(K_k)} (x-1)^{r(K_k)}
\\ \intertext{which, letting $\odd(k)=1$ if $k$ is odd and 0 otherwise}
   &= \sum_{k \leq n} \binom{n}{k} (y-1)^{\odd(k)} (x-1)^{k-\odd(k)}
\\ &= \sum_{k \operatorname{even}} \binom{n}{k} (x-1)^{k}
 + \frac{y-1}{x-1} \sum_{k \operatorname{odd}} \binom{n}{k} (x-1)^{k} .
\end{align*}
The even and odd sums are computable from the sum and difference of
$(z+1)^n = \sum \binom{n}{k} z^k$ and
$(-z+1)^n = \sum \binom{n}{k} z^k (-1)^k$;
substituting $z=x-1$ and simplifying gives the expression shown.

We also derive $q(K_{m,n})$ directly from \eqref{explicit}.
$K_{m,n}$ has $\binom{m}{i} \binom{n}{j}$ subgraphs $K_{i,j}$.
Each such subgraph's adjacency matrix has the form
\begin{align*}
A(K_{i,j}) &=
 \left( \begin{array}{cc} \zero & \one \\ \one & \zero \end{array} \right) ,
\intertext{whose rank is}
r(K_{i,j}) &=
  \begin{cases}
   0 & \textrm{if } i=0 \textrm{ or } j=0 \\
   2 & \textrm{if } i>0 \mbox{ and } j>0 .
  \end{cases}
\end{align*}
Then
\begin{align*}
q &= \sum (x-1)^r (y-1)^n
\\ &= \sum_{i=1}^m \sum_{j=1}^n
         (x-1)^2 (y-1)^{i+j-2} \binom{m}{i} \binom{n}{j}
      + \sum_{j=1}^n \binom{n}{j} + \sum_{i=1}^m \binom{m}{i} + 1 ,
\intertext{the four terms coming respectively from the cases where $i>0$ and
$j>0$; $i=0$ and $j>0$; $j=0$ and $i>0$; and $i=j=0$ (the null subgraph).
Expanding,}
 &= \frac{(x-1)^2}{(y-1)^2}
      \left( \sum_{i=1}^m (y-1)^i \binom{m}{i} \right)
      \left( \sum_{j=1}^n (y-1)^j \binom{n}{j} \right)
 \\ &\phantom{=}
    + \sum_{j=1}^n (y-1)^j \binom{n}{j}
    + \sum_{i=1}^m (y-1)^i \binom{m}{i} + 1 .
\end{align*}
The claim for $K_{m,n}$ follows immediately.

For $G= P_n$ with $n \geq 2$, we use the reduction~\eqref{red1}
with edge $ab$, where $b$ is a leaf.
Since $G-a$ is the disjoint union of
$P_{n-2}$ and $E_1$, $q(G-a) = y q(P_{n-2})$.
Since $\Gab=G$, $\Gab-b = P_{n-1}$ and $\Gab-a-b = P_{n-2}$.
The net result is
\begin{align*}
q(P_n) = (y+x^2-2x) q(P_{n-2}) + q(P_{n-1}) .
\end{align*}
Solving this recursion, with the boundary conditions
$q(P_0) = q(E_1) = y$ and $q(P_1) = q(K_2) = x^2-2x+2y$,
yields our formula for $q(P_n)$.
\end{proof}

\section{Further polynomials}
\label{otherpolys}

We observed in Section~\ref{intro} that the interlace polynomial's
expansion is similar to that of the Tutte polynomial,
with two significant differences:
the sum is over vertex rather than edge subsets,
and the rank and nullity are
those of an adjacency matrix rather than an incidence matrix.

This suggests a whole range of polynomials given by similar expansions,
with the sums taken variously over vertex or edge subsets,
with rank and nullity being those of an adjacency or an incidence matrix,
and with the field used perhaps being other than $\FF_2$.
Of course there is no obstacle to taking a ``master polynomial''
summing over both vertex and edge subsets,
incorporating terms for both the adjacency-matrix and the
incidence-matrix rank (using four variables instead of two),
and even introducing further variables to incorporate
ranks computed over other fields.
It would be interesting to determine which of these polynomials
satisfy reductions akin to that of Theorem~\ref{reduction},
and which ones have significance in combinatorics or other fields.

\section{Open problems}
\label{open}
As was the case with the single-variable interlace polynomial,
the two-variable interlace polynomial is quite new,
and there are more questions than answers.
Here we simply list a few of the obvious ones.

Is $q(G)$ reconstructible, \ie, given $q(G-a)$ for each vertex~$a$, can we
reconstruct $q(G)$?

What is the expectation of $q(G)$ for a random graph~$G$?

We conjectured in \cite{abs-jctb} that the \nullity polynomial's coefficient
sequence might be unimodal.
For the two-variable polynomial,
representing the coefficients of $q(G;-x,y)$ as an array
whose entry $(i,j)$ is the coefficient of $x^i y^j$,
the array's rows and columns are also unimodal,
for all (loopless) graphs through order~7.
(With the substitution of $-x$ for $x$,
it is clear from Corollary~\ref{implicit}
--- but not from \eqref{explicit} ---
that the coefficients are all non-negative.)
Unfortunately, there are six graphs of order 8 for which this is not the case.
Since however the vertex-rank polynomial and the vertex-nullity polynomial
for these six graphs
\emph{do} have unimodal coefficient sequences, it remains possible that
both of these single-variable polynomials' coefficient sequences
are always unimodal.
(As pointed out in \cite{abs-jctb} for the nullity polynomial, though,
there are reasons to be doubtful, including the falsification in 1993
of a similar conjecture for the Tutte polynomial~\cite{Schwarzler93}.)

But the most promising line of inquiry is
the proposal in the previous section,
to generalize from the interlace polynomial \eqref{explicit}
and the Tutte polynomial \eqref{Tutte}
to generate other interesting polynomials.

\section{Acknowledgments}
We are grateful to the two referees,
who brought more relevant literature to our attention
and whose comments resulted in structural and notational improvements.
We are also grateful to Hein van der Holst for providing us with
a manuscript copy of~\cite{AvdH} and for
helpful discussions on the interlace polynomials.

\bibliography{poly}

\end{document}